\documentclass[11pt,reqno,a4paper]{amsart}

\usepackage[T1]{fontenc}
\usepackage[utf8]{inputenc}
\usepackage{libertinus}
\usepackage{microtype}
\usepackage{geometry}
\usepackage{amsmath,amssymb,amsthm,amsfonts}
\usepackage{mathtools}
\usepackage{mathrsfs}
\usepackage{etoolbox}
\usepackage{aliascnt}
\numberwithin{equation}{section}

\usepackage{xcolor}
\usepackage[
  colorlinks=true,
  linkcolor=blue,
  citecolor=blue,
  urlcolor=blue,
  pagebackref=true
]{hyperref}
\renewcommand*{\backrefalt}[4]{%
  \ifcase #1\relax
  \or {\footnotesize Cited on page #2.}%
  \else {\footnotesize Cited on pages #2.}%
  \fi
}
\usepackage[nameinlink,noabbrev]{cleveref}
\crefname{section}{Section}{Sections}
\Crefname{section}{Section}{Sections}
\crefname{subsection}{Subsection}{Subsections}
\Crefname{subsection}{Subsection}{Subsections}
\crefname{theorem}{Theorem}{Theorems}
\Crefname{theorem}{Theorem}{Theorems}
\crefname{proposition}{Proposition}{Propositions}
\Crefname{proposition}{Proposition}{Propositions}
\crefname{lemma}{Lemma}{Lemmas}
\Crefname{lemma}{Lemma}{Lemmas}
\crefname{corollary}{Corollary}{Corollaries}
\Crefname{corollary}{Corollary}{Corollaries}
\crefname{maintheorem}{Theorem}{Theorems}
\Crefname{maintheorem}{Theorem}{Theorems}
\crefname{pointwiseremark}{Remark}{Remarks}
\Crefname{pointwiseremark}{Remark}{Remarks}

\usepackage{titlesec}
\titleformat{\section}
  {\normalfont\large\bfseries\centering}
  {\thesection.}{0.75em}{}
\titlespacing*{\section}{0pt}{3.2ex plus 1ex minus .2ex}{1.6ex plus .3ex}
\titleformat{\subsection}
  {\normalfont\normalsize\bfseries}
  {\thesubsection.}{0.75em}{}
\titlespacing*{\subsection}{0pt}{2.4ex plus .8ex minus .2ex}{1.0ex plus .2ex}
\newtheoremstyle{compact}
  {6pt}{1pt}{\itshape}{}{\bfseries}{.}{0.5em}{}
\theoremstyle{compact}
\newtheorem{theorem}{Theorem}[section]
\newaliascnt{proposition}{theorem}
\newtheorem{proposition}[proposition]{Proposition}
\aliascntresetthe{proposition}
\newaliascnt{lemma}{theorem}

\aliascntresetthe{lemma}
\newaliascnt{corollary}{theorem}
\newtheorem{corollary}[corollary]{Corollary}
\aliascntresetthe{corollary}
\newtheorem{maintheorem}{Theorem}

\theoremstyle{definition}
\newtheorem{remark}[theorem]{Remark}
\newtheorem{pointwiseremark}{Remark}[section]

\AtBeginEnvironment{proof}{\vspace{6pt}}

\makeatletter
\AtBeginEnvironment{thebibliography}{%
  \def\@mklab#1{#1\hfil}%
}
\makeatother

\newcommand{\R}{\mathbb R}
\newcommand{\cS}{\mathcal S}
\newcommand{\cP}{\mathcal P}
\newcommand{\supp}{\operatorname{supp}}
\newcommand{\dimF}{\dim_{\mathrm F}}
\newcommand{\dimH}{\dim_{\mathrm H}}

\begin{document}

\title{Ergodic theorems for homogeneous dilations}

\author{Michael Bj\"orklund}
\address{Department of Mathematics, Chalmers University of Technology and University of Gothenburg, Gothenburg, Sweden}
\email{micbjo@chalmers.se}

\begin{abstract}
Let $\nu$ be a compactly supported probability measure on $\R^d$ with Fourier
dimension greater than one.  We prove an almost-everywhere ergodic theorem for
the corresponding dilated averages under ergodic measure-preserving
$\R^d$-actions.  If $a=\dimF\nu$, the result holds on $L^p$ for
$1+1/a<p<\infty$.  The proof combines Rubio de Francia's Euclidean maximal
theorem with a transference argument and a dense class obtained by spectral
localization away from the origin.  Applications include curved
hypersurfaces and averaging measures carried by Salem sets, including
Brownian images and explicit Diophantine constructions.
\end{abstract}

\maketitle

\section{Introduction}
\label{sec:introduction}

Let $(X,\mathcal B,\mu)$ be a standard probability space and let
$T=(T_t)_{t\in\R^d}$ be a jointly Borel, ergodic, measure-preserving action of
$\R^d$.  Write $U_t f:=f\circ T_t$.  Throughout, $1\leq p<\infty$.  For a
compactly supported probability measure $\nu$ on $\R^d$ and $\lambda>0$, set
\begin{equation}\label{eq:averages}
  A_\lambda f:=\int_{\R^d} U_{\lambda t}f\,d\nu(t)
  \qquad\text{in }L^p(X).
\end{equation}
For each fixed $\lambda>0$, the operator in \eqref{eq:averages} has the
expected almost-everywhere pointwise interpretation.  The only minor issue is
to choose these representatives consistently when $\lambda$ varies; see
\cref{rem:varying-scales} below.

We use the Fourier transform
\[
  \widehat\nu(\xi):=\int_{\R^d}e^{-2\pi i\langle \xi,t\rangle}\,d\nu(t)
\]
and define the Fourier dimension of $\nu$ by
\begin{equation}\label{eq:fourier-dimension}
  \dimF\nu:=\sup\bigl\{0\leq s\leq d:
  |\widehat\nu(\xi)|\lesssim_s (1+|\xi|)^{-s/2}\bigr\}.
\end{equation}

Our main result is the following.

\begin{maintheorem}\label{thm:main}
Let $\nu$ be a compactly supported probability measure on $\R^d$ and suppose
that
\[
  a:=\dimF\nu>1.
\]
Then, for every ergodic measure-preserving action $T$ as above, every
\[
  1+\frac1a<p<\infty,
\]
and every $f\in L^p(X)$,
\[
  \lim_{\lambda\to\infty} A_\lambda f(x)=\int_X f\,d\mu
\]
for almost every $x\in X$.
\end{maintheorem}

Notice that $a>1$ forces $d\geq2$.

The basic absolutely continuous example is normalized Lebesgue measure on the
unit ball.  In that case the operators $A_\lambda$ are precisely Wiener's ball
averages, and Wiener's ergodic theorem gives almost-everywhere convergence for
all $f\in L^1(X)$ \cite{Wie39}.  Thus \cref{thm:main} is not meant to improve
the classical theorem in the absolutely continuous case; it extends the same
averaging picture to singular measures with sufficiently large Fourier
dimension.

The Euclidean maximal estimate is due to Rubio de Francia \cite{Rub86}.  If
$\nu$ is compactly supported and
\begin{equation}\label{eq:intro-decay}
  |\widehat\nu(\xi)|\lesssim (1+|\xi|)^{-b/2},
  \qquad b>1,
\end{equation}
then the maximal function, initially defined for $F\in\cS(\R^d)$ by
\[
  \mathcal M_\nu F:=\sup_{\lambda>0}
  \left|\int_{\R^d}F(\,\cdot+\lambda t)\,d\nu(t)\right|,
\]
satisfies $\|\mathcal M_\nu F\|_p\lesssim_{p,\nu}\|F\|_p$ for
$1+1/b<p<\infty$.  We give the finite-scale
transference argument in full in \cref{sec:maximal}; after passage to a
countable dense set of scales it gives the ergodic maximal inequality needed
below.

The pointwise input is proved independently of the maximal theorem.  We use
smoothings whose spectral support is a compact subset of
$\R^d\setminus\{0\}$.  These functions form a dense subspace of
$L^p_0(X)$.  If $f$ belongs to this subspace and \eqref{eq:intro-decay} holds,
the spectral theorem gives
\[
  \|A_\lambda f\|_2^2\lesssim_{\nu,f} \lambda^{-b},
  \qquad \lambda\geq1.
\]
Thus $\lambda\mapsto A_\lambda f(x)$ belongs to $L^2([1,\infty))$ for almost
every $x$.  Compact support of $\nu$ makes the same map uniformly Lipschitz.
An $L^2$ Lipschitz function on $[1,\infty)$ tends to zero.  At this stage $b>1$ is exactly the integrability threshold for
$\lambda^{-b}$; it is also the threshold required by Rubio de Francia's
maximal theorem.

The passage from a decay exponent $b$ to Fourier dimension requires no
endpoint assumption.  If $a=\dimF\nu>1$ and $p>1+1/a$, choose
$1<b<a$ sufficiently close to $a$ that $p>1+1/b$.  The definition of Fourier
dimension then gives \eqref{eq:intro-decay} with this $b$.

The restriction $a>1$ is a restriction of this argument.  In particular,
normalized arc-length measure on $S^1\subset\R^2$ has Fourier dimension one,
so \cref{thm:main} does not contain the planar circle ergodic theorem.  That
endpoint rests on Bourgain's circular maximal theorem \cite{Bou86}; see
Lacey \cite{Lac95} for the corresponding ergodic theorem.

\subsection*{Erratum on a previous version}
An earlier version of this work (arXiv:0902.0680) was withdrawn in June 2011
after the author found a gap in the dense-subspace argument for pointwise
convergence.  The present paper should be read as a correction to the
Fourier-dimension theorem in that version and gives a new proof of that result.
The broader statements there on Rajchman measures under homogeneous dilations
and on polynomial curves are not addressed here.

\section{A dense class with pointwise convergence}
\label{sec:dense}

\begin{pointwiseremark}[Pointwise interpretation and varying scales]
\label{rem:varying-scales}
For a fixed $\lambda>0$ and a Borel representative of $f\in L^p(X)$,
Tonelli and invariance give
\[
  \int_X\int_{\R^d}|f(T_{\lambda t}x)|\,d\nu(t)\,d\mu(x)=\|f\|_1.
\]
Thus the literal average is finite almost everywhere.  Minkowski's inequality
and duality identify its $L^p$-class with the Bochner integral in
\eqref{eq:averages}; the same argument applied to the indicator of a null set
shows that this class is independent of the representative of $f$.

The exceptional set in this fixed-scale identification may depend on
$\lambda$, so one cannot intersect these sets over all $\lambda>0$.  When a
pointwise statement below involves a varying scale, $A_\lambda f(x)$ therefore
denotes the simultaneous measurable realization constructed in
\cref{sec:proof-main}.  For each fixed $\lambda$ it represents
\eqref{eq:averages}.  The construction is continuous in $\lambda$ almost
everywhere; that continuity is used only to pass from rational to arbitrary
scales.
\end{pointwiseremark}

We now recall the standard spectral notation.  On $L^2(X)$ let $E$ denote the
projection-valued spectral measure of the Koopman representation, normalized
so that
\[
  U_t=\int_{\R^d}e^{-2\pi i\langle t,\xi\rangle}\,dE(\xi).
\]
We use without further comment the spectral theorem and functional calculus
for unitary representations of abelian groups; see
\cite[Corollary~2.13 and \S2.6]{EW25}.  For $\varphi\in L^1(\R^d)$ set
\[
  U(\varphi)g:=\int_{\R^d}\varphi(s)U_sg\,ds,
\]
so that
\begin{equation}\label{eq:spectral-multipliers}
  U(\varphi)=\int_{\R^d}\widehat\varphi(\xi)\,dE(\xi),
  \qquad
  A_\lambda=\int_{\R^d}\widehat\nu(\lambda\xi)\,dE(\xi).
\end{equation}
For $g\in L^2(X)$ write
$\sigma_g(B):=\langle E(B)g,g\rangle$.  Then
\[
  d\sigma_{U(\varphi)g}=|\widehat\varphi|^2\,d\sigma_g;
\]
compare \cite[Exercise~2.60]{EW25}.  If $g$ is bounded and Borel, the formula
\[
  x\longmapsto\int_{\R^d}\varphi(s)g(T_sx)\,ds
\]
gives a bounded Borel representative of $U(\varphi)g$, by Fubini.
Let $L^p_0(X)$ denote the mean-zero subspace of $L^p(X)$.

\begin{proposition}\label{prop:dense}
Let $T$ be ergodic and fix $1\leq p<\infty$.  The linear span of the functions
\[
  U(\varphi)g,
  \qquad
  g\in L^\infty(X),\quad
  \varphi\in\cS(\R^d),\quad
  \supp\widehat\varphi\Subset\R^d\setminus\{0\},
\]
is dense in $L^p_0(X)$.  Every function in this span has mean zero.
\end{proposition}

\begin{proof}
Since $\widehat\varphi(0)=0$, invariance of $\mu$ gives
$\int_XU(\varphi)g\,d\mu=0$.  It remains to prove density.  Bounded mean-zero
functions are dense in $L^p_0(X)$: truncate an arbitrary $h\in L^p_0(X)$ and
subtract the mean of the truncation.  It is therefore enough to fix
$g\in L^\infty_0(X)$.  Choose
$\eta\in\cS(\R^d)$ such that $\widehat\eta\in C_c^\infty(\R^d)$,
$\widehat\eta=1$ on a neighbourhood of the origin, and
$\widehat\eta(0)=1$.  For $r>0$, put
$\eta_r(s):=r^{-d}\eta(s/r)$.

The spectral identities in \eqref{eq:spectral-multipliers} give
\[
  \|U(\eta_r)g-g\|_2^2
  =\int_{\R^d}|\widehat\eta(r\xi)-1|^2\,d\sigma_g(\xi)
  \longrightarrow0
\]
as $r\downarrow0$.  By the spectral theorem, $E(\{0\})$ is the orthogonal
projection onto the subspace corresponding to the trivial character, hence
onto the $U$-invariant vectors; see \cite[\S2.6.1]{EW25}.  Ergodicity says that
this invariant subspace consists of the constants (cf.\ \cite[\S6.2.2]{EW25}).
Since $g$ has mean zero, $\sigma_g(\{0\})=0$, and therefore
\[
  \|U(\eta_R)g\|_2^2
  =\int_{\R^d}|\widehat\eta(R\xi)|^2\,d\sigma_g(\xi)
  \longrightarrow0
\]
as $R\to\infty$.  The operators $U(\eta_r)$ are uniformly bounded on
$L^\infty$, so these two $L^2$ limits imply the corresponding $L^p$ limits for
every finite $p$ (by the probability-space inclusion for $p\leq2$ and
interpolation with $L^\infty$ for $p>2$).

For $R$ large and $r$ small, $U(\eta_r-\eta_R)g$ is therefore arbitrarily
close to $g$ in $L^p(X)$.  Moreover,
$\widehat\eta(r\xi)-\widehat\eta(R\xi)$ has compact support and vanishes on a
neighbourhood of the origin.  This proves the assertion.
\end{proof}

\begin{proposition}\label{prop:dense-convergence}
Let $T$ be a jointly Borel measure-preserving action of $\R^d$ on a standard
probability space.  Suppose that $\nu$ is compactly supported and that, for some
$b>1$,
\begin{equation}\label{eq:decay-b}
  |\widehat\nu(\xi)|\leq C(1+|\xi|)^{-b/2}.
\end{equation}
Then
\[
  A_\lambda f(x)\longrightarrow0
\]
for almost every $x$ and every $f$ in the dense class of
\cref{prop:dense}.
\end{proposition}

\begin{proof}
It is enough to consider a generator
$f=U(\varphi)g$, where $g\in L^\infty(X)$ and
$\supp\widehat\varphi\Subset\R^d\setminus\{0\}$.  We use the bounded Borel
representative
\[
  f(x)=\int_{\R^d}\varphi(s)g(T_sx)\,ds.
\]
For this representative the averages form a jointly Borel family, obtained by
writing the literal average as the absolutely convergent double integral
\[
  A_\lambda f(x)
  =\int_{\R^d}\int_{\R^d}
    \varphi(s)g(T_{s+\lambda t}x)\,ds\,d\nu(t).
\]
For each fixed $\lambda$ this represents the operator in
\eqref{eq:averages}, by \cref{rem:varying-scales}.

\smallskip
\noindent\emph{1. Square integrability in the scale parameter.}
Let
\[
  \delta:=\operatorname{dist}
  \bigl(0,\supp\widehat\varphi\bigr)>0.
\]
The first identity below is simply the $L^2$ norm formula for the spectral
multiplier $A_\lambda$.  Indeed, by \eqref{eq:spectral-multipliers},
$A_\lambda$ has multiplier
$m_\lambda(\xi)=\widehat\nu(\lambda\xi)$; hence the functional calculus and
the definition $\sigma_f(B)=\langle E(B)f,f\rangle$ give
\[
  \|A_\lambda f\|_2^2
  =\bigl\langle |m_\lambda|^2(E)f,f\bigr\rangle
  =\int_{\R^d}|\widehat\nu(\lambda\xi)|^2\,d\sigma_f(\xi).
\]
On the other hand,
$d\sigma_f=|\widehat\varphi|^2\,d\sigma_g$, so $\sigma_f$ is supported on
$\supp\widehat\varphi\subset\{\xi:|\xi|\geq\delta\}$.  Thus every frequency
seen by $f$ is separated from the origin, and the Fourier decay
\eqref{eq:decay-b} is uniform on the support of $\sigma_f$.  Since
$\sigma_f(\R^d)=\|f\|_2^2$, we obtain
\begin{equation}\label{eq:l2-scale}
\begin{aligned}
  \|A_\lambda f\|_2^2
  &=\int_{\R^d}|\widehat\nu(\lambda\xi)|^2\,d\sigma_f(\xi)\\
  &\leq C^2(1+\lambda\delta)^{-b}\,\sigma_f(\R^d)\\
  &=C^2(1+\lambda\delta)^{-b}\,\|f\|_2^2
  \lesssim_{\nu,\varphi}\lambda^{-b}\|f\|_2^2,
  \qquad \lambda\geq1.
\end{aligned}
\end{equation}
This is precisely why the dense class is localized away from frequency zero:
under dilation, the whole spectral support is pushed into the region where
$\widehat\nu$ has quantitative decay.  Since $b>1$, integrating
\eqref{eq:l2-scale} in $\lambda$ and applying Fubini gives
\begin{equation}\label{eq:pointwise-l2-scale}
  \int_1^\infty |A_\lambda f(x)|^2\,d\lambda<\infty
\end{equation}
for almost every $x$.  By itself this does not imply convergence as
$\lambda\to\infty$: an $L^2$ function may have arbitrarily thin spikes.  The
next estimate rules out exactly that possibility.

\smallskip
\noindent\emph{2. Uniform Lipschitz regularity in the scale parameter.}
Changing variables in the double integral above gives
\[
  A_\lambda f(x)
  =\int_{\R^d}g(T_sx)K_\lambda(s)\,ds,
  \qquad
  K_\lambda(s):=\int_{\R^d}\varphi(s-\lambda t)\,d\nu(t).
\]
By the mean value inequality,
\[
  \|\varphi(\,\cdot-\lambda t)-\varphi(\,\cdot-\rho t)\|_1
  \leq |\lambda-\rho|\,|t|\,\|\nabla\varphi\|_1.
\]
After integration against $\nu$ this yields
\begin{equation}\label{eq:lipschitz-kernel}
  \|K_\lambda-K_\rho\|_1
  \leq |\lambda-\rho|\,\|\nabla\varphi\|_1
  \int_{\R^d}|t|\,d\nu(t).
\end{equation}
(The last integral is finite because $\nu$ is compactly supported.)  It follows
that, for every $x$, the map $\lambda\mapsto A_\lambda f(x)$ is Lipschitz with
a constant independent of $x$.

\smallskip
\noindent\emph{3. Conclusion.}
A Lipschitz function $F\in L^2([1,\infty))$ must tend to zero.  Indeed, if
$F$ is $K$-Lipschitz and $K>0$, then
\[
  |F(\ell)|^3
  \leq 8K\int_\ell^\infty |F(t)|^2\,dt,
  \qquad \ell\geq1,
\]
because $|F(t)|\geq|F(\ell)|/2$ on
$[\ell,\ell+|F(\ell)|/(2K)]$.  The right-hand side tends to zero; the case
$K=0$ is immediate.  Combining this observation with
\eqref{eq:pointwise-l2-scale} proves $A_\lambda f(x)\to0$ almost everywhere.
Finite linear combinations are handled identically.
\end{proof}

\section{The maximal theorem}
\label{sec:maximal}

For $F\in\cS(\R^d)$ define
\[
  \mathcal M_\nu F(y):=
  \sup_{\lambda>0}\left|\int_{\R^d}F(y+\lambda t)\,d\nu(t)\right|.
\]
Rubio de Francia's theorem is the only Euclidean maximal result used below.

The following is \cite[Theorem A]{Rub86} applied, with our sign convention,
to the reflected measure.

\begin{theorem}[Rubio de Francia]\label{thm:rubio}
Suppose that $\nu$ is a finite compactly supported Borel measure on $\R^d$
satisfying \eqref{eq:decay-b} for some $b>1$.  Then, for every
$p>1+1/b$,
\[
  \|\mathcal M_\nu F\|_{L^p(\R^d)}
  \leq C_{p,\nu}\|F\|_{L^p(\R^d)},
  \qquad F\in\cS(\R^d).
\]
\end{theorem}

More precisely, Rubio de Francia writes the averages with
$F(y-\lambda t)$.  For $\check\nu(E):=\nu(-E)$ one has
$|\widehat{\check\nu}|=|\widehat\nu|$, which yields the displayed
$F(y+\lambda t)$ formulation.

Write $\mathbb Q_+:=\mathbb Q\cap(0,\infty)$.  The Euclidean maximal
operator above uses all real scales; on $X$ we first use rational scales so
that the supremum is representative-independent.  For $f\in L^p(X)$ define
the separable maximal function
\[
  \mathcal M_\nu^T f:=\sup_{\lambda\in\mathbb Q_+}|A_\lambda f|.
\]
Choose Borel representatives of the countable family
$(A_\lambda f)_{\lambda\in\mathbb Q_+}$.  Changing these choices changes the
supremum only on a countable union of null sets, so
$\mathcal M_\nu^T f$ is a well-defined measurable $L^0$-class.

\begin{theorem}\label{thm:ergodic-maximal}
Suppose that $\nu$ is a compactly supported probability measure satisfying
\eqref{eq:decay-b} for some $b>1$.  Let $T$ be a jointly Borel,
measure-preserving action of $\R^d$ on a standard probability space.  If
$p>1+1/b$, then
\begin{equation}\label{eq:ergodic-maximal}
  \|\mathcal M_\nu^T f\|_{L^p(X)}
  \leq C_{p,\nu}\|f\|_{L^p(X)}.
\end{equation}
\end{theorem}

The constant in \eqref{eq:ergodic-maximal} may be taken to be the Rubio de
Francia constant; in particular it is independent of the finite scale sets
and truncation radii used below.

\begin{proof}
By \cref{thm:rubio}, let $C_{p,\nu}$ be a bound for $\mathcal M_\nu$ on
$L^p(\R^d)$.  Choose $M>0$ with $\supp\nu\subset B_M$.  It is enough to prove
the same estimate, with the same constant, for the maximum over an arbitrary
finite set $\Lambda\subset\mathbb Q_+$.  Put $L:=\max\Lambda$.

Fix $R>0$ and a Borel representative of $f$.  For $x\in X$ define
\[
  F_x(s):=f(T_sx)\,\mathbf 1_{B_{R+LM}}(s).
\]
By Fubini, $F_x\in L^p(\R^d)$ for almost every $x$.  We first record the
fixed-scale interpretation of convolution for arbitrary $L^p(\R^d)$ data.
For a Borel representative of $F\in L^p(\R^d)$ and fixed $\lambda>0$, the
literal integral
\[
  C_\lambda F(s):=\int_{\R^d}F(s+\lambda t)\,d\nu(t)
\]
exists for almost every $s$.  Indeed, for every ball $B$,
\[
  \int_B\int_{\R^d}|F(s+\lambda t)|\,d\nu(t)\,ds
  \leq |B|^{1-1/p}\|F\|_p<\infty
\]
(with the evident $p=1$ interpretation).  The same calculation with the
indicator of a null set shows that this a.e.-defined function is independent
of the chosen representative.  Minkowski's inequality gives
$\|C_\lambda F\|_p\leq\|F\|_p$.

For a finite family $\Lambda$, define
\[
  \mathcal M_{\nu,\Lambda}F
  :=\max_{\lambda\in\Lambda}|C_\lambda F|.
\]
For Schwartz functions,
$\mathcal M_{\nu,\Lambda}F\leq\mathcal M_\nu F$.  If
$F_k\in\cS(\R^d)$ and $F_k\to F$ in $L^p$, Minkowski gives
$C_\lambda F_k\to C_\lambda F$ in $L^p$ for every $\lambda\in\Lambda$.
After passing to one subsequence, these convergences hold almost everywhere
simultaneously for all $\lambda\in\Lambda$.  Fatou's lemma and
\cref{thm:rubio} therefore yield
\[
  \|\mathcal M_{\nu,\Lambda}F\|_p
  \leq \liminf_k\|\mathcal M_{\nu,\Lambda}F_k\|_p
  \leq \liminf_k\|\mathcal M_\nu F_k\|_p
  \leq C_{p,\nu}\|F\|_p.
\]

For almost every $x$, apply this to $F_x$.  If $s\in B_R$,
$\lambda\in\Lambda$, and $t\in\supp\nu$, then
$s+\lambda t\in B_{R+LM}$.  For each fixed $\lambda\in\Lambda$, the action identity and the
fixed-scale identification after \eqref{eq:averages} give this equality away
from a $\mu$-null exceptional set in the orbit variable.  Invariance of
$\mu$ and Fubini transfer that exceptional set to a null subset of
$B_R\times X$.  Hence, for almost every $(s,x)\in B_R\times X$,
\[
  A_\lambda f(T_sx)=C_\lambda F_x(s),
  \qquad \lambda\in\Lambda.
\]
Since $\Lambda$ is finite, the identities hold simultaneously outside one
null set.  Therefore
\[
  \max_{\lambda\in\Lambda}|A_\lambda f(T_sx)|
  = \mathcal M_{\nu,\Lambda}F_x(s)
\]
for almost every $(s,x)\in B_R\times X$.  Using invariance of $\mu$ and
Fubini's theorem,
we obtain
\begin{align*}
  |B_R|\left\|\max_{\lambda\in\Lambda}|A_\lambda f|\right\|_p^p
  &=\int_{B_R}\int_X
    \max_{\lambda\in\Lambda}|A_\lambda f(T_sx)|^p\,d\mu(x)\,ds\\
  &\leq \int_X\int_{\R^d}|\mathcal M_{\nu,\Lambda} F_x(s)|^p\,ds\,d\mu(x)\\
  &\leq C_{p,\nu}^p
    \int_X\int_{B_{R+LM}}|f(T_sx)|^p\,ds\,d\mu(x)\\
  &=C_{p,\nu}^p |B_{R+LM}|\,\|f\|_p^p.
\end{align*}
Letting $R\to\infty$ gives
\[
  \left\|\max_{\lambda\in\Lambda}|A_\lambda f|\right\|_p
  \leq C_{p,\nu}\|f\|_p.
\]
Exhaust $\mathbb Q_+$ by finite sets and apply monotone convergence to obtain
\eqref{eq:ergodic-maximal}.
\end{proof}

\section{Proof of the pointwise theorem}
\label{sec:proof-main}

\begin{proof}[Proof of \cref{thm:main}]
Fix $1+1/a<p<\infty$.  Choose $1<b<a$ such that
\[
  p>1+\frac1b.
\]
By the definition of Fourier dimension, \eqref{eq:decay-b} holds for this
$b$.  Hence \cref{thm:ergodic-maximal} gives
\begin{equation}\label{eq:maximal-rational}
  \|\mathcal M_\nu^T h\|_p\lesssim_{p,\nu}\|h\|_p.
\end{equation}

It remains to work in $L^p_0(X)$.  Fix $f\in L^p_0(X)$ and a Borel
representative of $f$.  Choose $g_n$ in the dense class from
\cref{prop:dense}, represented by their defining smoothing integrals, such
that $g_n\to f$ in $L^p$.  Passing to a subsequence, \eqref{eq:maximal-rational}
allows us to arrange
\[
  \sum_{n=1}^\infty
  \|\mathcal M_\nu^T(f-g_n)\|_p^p<\infty.
\]
For the chosen Borel representatives, define the rational-scale maximal
function by the literal integrals
\[
  H_n(x):=\sup_{q\in\mathbb Q_+}
  \left|\int_{\R^d}(f-g_n)(T_{qt}x)\,d\nu(t)\right|,
\]
with the value $+\infty$ if one of the displayed integrals is not absolutely
convergent.  On the common conull set used below all these integrals are
finite.  For each $n$, this is a representative of
$\mathcal M_\nu^T(f-g_n)$: each fixed rational-scale integral represents the
corresponding $L^p$ operator, and there are only countably many scales.
Consequently $H_n(x)\to0$ for almost every $x$.

Intersect this conull set with the conull sets on which
\cref{prop:dense-convergence} holds for every $g_n$, and with the conull set on
which all rational-scale integrals above are finite.  Fix $x$ in the resulting
set.  Each map $\lambda\mapsto A_\lambda g_n(x)$ is continuous on
$(0,\infty)$ and tends to zero as $\lambda\to\infty$.  For $m,n\geq1$ and
$q\in\mathbb Q_+$, pointwise linearity of the chosen integral representatives
gives
\[
  |A_q(g_n-g_m)(x)|\leq H_n(x)+H_m(x).
\]
The left-hand side is continuous as a function of the scale.  Since
$\mathbb Q_+$ is dense in $(0,\infty)$,
\[
  \sup_{\lambda>0}|A_\lambda(g_n-g_m)(x)|
  \leq H_n(x)+H_m(x).
\]
Thus $(A_\lambda g_n(x))_n$ is uniformly Cauchy on $(0,\infty)$.  Let
$F_x(\lambda)$ denote its uniform limit.  Then $F_x$ is continuous and, on
letting $m\to\infty$ in the preceding inequality,
\begin{equation}\label{eq:uniform-version-error}
  \sup_{\lambda>0}|F_x(\lambda)-A_\lambda g_n(x)|\leq H_n(x).
\end{equation}
Since $H_n(x)\to0$ and $A_\lambda g_n(x)\to0$ as $\lambda\to\infty$ for each
$n$, \eqref{eq:uniform-version-error} implies
$F_x(\lambda)\to0$ as $\lambda\to\infty$.

Take the conull set above to be a measurable set $\Omega$ and define
\[
  F(\lambda,x):=\lim_{n\to\infty}
  \mathbf 1_{\Omega}(x)A_\lambda g_n(x).
\]
The limit exists for every $(\lambda,x)$: on $\Omega$ it is the uniform limit
constructed above, and off $\Omega$ it is zero.  Each term is jointly Borel by
\cref{prop:dense-convergence}, so $F$ is jointly measurable.  It remains to
identify this family with the original
$L^p$ averages.  For each fixed $\lambda>0$, the contraction property gives
$A_\lambda g_n\to A_\lambda f$ in $L^p(X)$, hence in measure.  On the conull
set fixed above, $A_\lambda g_n(x)\to F_x(\lambda)$.  Uniqueness of limits in
measure therefore shows that $x\mapsto F_x(\lambda)$ represents
$A_\lambda f$ for every fixed $\lambda$.  We take this jointly measurable
family $F(\lambda,x)$ as the simultaneous realization specified in the
convention of \cref{sec:introduction}, and hence write it simply as
$A_\lambda f(x)$.  This proves the theorem for mean-zero functions.
Finally, $A_\lambda 1=1$.  Applying the mean-zero result to
$f-\int_Xf\,d\mu$ and then adding the constant
$\int_Xf\,d\mu$ completes the proof.
\end{proof}

\section{Examples}
\label{sec:examples}

The theorem is useful when Fourier decay is obtained from geometry or from a
fractal construction.  We record several instances.

\subsection{Curved hypersurfaces}
Let $S\subset\R^d$ be a smooth hypersurface with non-vanishing Gaussian
curvature, and let $\nu$ be a compactly supported smooth probability measure
on $S$.  Stationary phase gives
\[
  |\widehat\nu(\xi)|\lesssim (1+|\xi|)^{-(d-1)/2};
\]
see \cite[Chapter VIII]{Ste93}.  Since the Fourier dimension of a measure is bounded by the Hausdorff
dimension of its support \cite[\S3.6 and Theorem~3.10]{Mat15}, one has
$\dimF\nu=d-1$.  For $d\geq3$,
\cref{thm:main} applies for
\[
  \frac d{d-1}<p<\infty.
\]
For normalized surface measure on the sphere this is the same $L^p$ range as
Jones' ergodic theorem \cite{Jon93}.  In dimension two the Fourier dimension
is one, and the argument above does not reach the circle theorem
\cite{Lac95}.

\subsection{Salem sets}
For a compact set $E\subset\R^d$, define
\[
  \dimF E:=\sup\{\dimF\nu:\nu\in\cP(E)\},
\]
where $\cP(E)$ is the set of Borel probability measures supported on $E$.
One always has $\dimF E\leq\dimH E$; a set for which equality holds is called
a \emph{Salem set}.  See \cite[\S3.6 and Theorem~3.10]{Mat15}.

\begin{corollary}\label{cor:salem-set}
Let $E\subset\R^d$ be a compact Salem set of Hausdorff dimension $s>1$.
For every
\[
  1+\frac1s<p<\infty
\]
there is a probability measure $\nu$ supported on $E$ such that, for every
jointly Borel ergodic measure-preserving $\R^d$-action and every
$f\in L^p(X)$, the family
\[
  A_\lambda f:=\int_E U_{\lambda t}f\,d\nu(t)
  \qquad\text{in }L^p(X)
\]
satisfies
\[
  A_\lambda f(x)\longrightarrow \int_Xf\,d\mu
\]
for almost every $x$ as $\lambda\to\infty$, with the convention of
\cref{sec:introduction}.
\end{corollary}

\begin{proof}
Choose $b<s$ with $b>1$ and $p>1+1/b$.  Since $\dimF E=s$, there is
$\nu\in\cP(E)$ with $\dimF\nu>b$.  Apply \cref{thm:main}.
\end{proof}

The measure in \cref{cor:salem-set} may depend on $p$; the definition of a
Salem set does not require a single measure to attain the endpoint exponent
$s$.

\subsection{Explicit Diophantine Salem sets}
Fraser and Hambrook \cite[Theorems 1.1--1.2]{FH23} construct explicit Salem sets in
$\R^d$ of every prescribed dimension $0<s<d$, using Diophantine approximation
in number fields.  Their construction produces compactly supported
probability measures of Fourier dimension $s$ on these sets.  The endpoint
Fourier estimate has a subpolynomial loss, which does not affect the Fourier
dimension or the application of \cref{thm:main}.  Thus, for every
$d\geq2$, every $1<s<d$, and every $1+1/s<p<\infty$, one obtains explicit singular
averaging measures to which \cref{thm:main} applies.

\subsection{Fractional Brownian and Gaussian images}
Let $B^\alpha:\R^N\to\R^d$ be fractional Brownian motion of Hurst index
$0<\alpha<1$, and fix a compact set $E\subset\R^N$ with
$s:=\dimH E\leq\alpha d$.  By \cite[Corollary 2.10]{SX06}, for this fixed
$E$ the image $B^\alpha(E)$ is almost surely a Salem set with
\[
  \dimF B^\alpha(E)=\dimH B^\alpha(E)=\frac{s}{\alpha}.
\]
Consequently, if $s>\alpha$, then for almost every sample path and for each
$1+\alpha/s<p<\infty$, \cref{cor:salem-set} gives a probability measure supported on
$B^\alpha(E)$ to which \cref{thm:main} applies.  For ordinary Brownian motion
($N=1$, $\alpha=1/2$) in $\R^d$ with $d\geq2$, a compact
$E\subset[0,1]$ of dimension $s>1/2$ therefore gives the range
$1+1/(2s)<p<\infty$; here $s\leq1\leq d/2$, so the condition
$s\leq\alpha d$ is automatic.  The same paper treats substantially more
general Gaussian random fields; see \cite[Theorem 2.12]{SX06}.

\subsection{Fractal percolation and multiplicative cascades}
Consider dyadic fractal percolation in $[0,1]^d$ with retention probability
$q>2^{-d}$, so that non-extinction has positive probability.  On
non-extinction its limit set and natural measure have almost sure Hausdorff
dimension
\[
  s=d+\log_2 q.
\]
If $0\leq s\leq2$, Shmerkin and Suomala
\cite[Theorem 14.1]{SS18} show that the natural measure is almost surely a
Salem measure on the non-extinction event.  After normalizing its total mass,
we therefore obtain, whenever $1<s\leq2$, a compactly supported probability
measure of Fourier dimension $s$ to which \cref{thm:main} applies for
$1+1/s<p<\infty$.  For $q<1$, this is also the Bernoulli special case of the cascade theorem
below.

Chen, Li and Suomala \cite[Theorem 1.2]{CLS25} treat general dyadic
Mandelbrot multiplicative cascades.  Let $W\geq0$ satisfy
$\mathbb E W=1$, $\mathbb E(W\log_2 W)<d$, assume that $W$ is not almost
surely equal to one, and suppose that $W$ is sub-exponential in the sense that
there is $c>0$ such that
$\mathbb P(W>t)\leq 2e^{-ct}$ for all $t>0$.  If $\mu$ is the corresponding
cascade measure, define its
correlation dimension by
\[
  \dim_2\mu
  :=\liminf_{n\to\infty}
  \frac{-\log_2\!\left(\sum_{Q\in\mathcal D_n}\mu(Q)^2\right)}{n},
\]
where $\mathcal D_n$ is the family of dyadic cubes of side length $2^{-n}$.
Then, almost surely on non-extinction,
\[
  \dimF\mu=\min\{2,\dim_2\mu\}.
\]
Thus, whenever $\dim_2\mu>1$, normalizing $\mu$ to total mass one gives,
almost surely on non-extinction, a compactly supported multifractal averaging
measure covered by \cref{thm:main}.  The excluded degenerate case
$W\equiv1$ gives normalized Lebesgue measure on $[0,1]^d$ and can be treated
separately.

\begin{remark}\label{rem:hypotheses}
The proof uses the two quantitative hypotheses in different places.  Fourier
decay with exponent $b>1$ gives the integrable estimate
\eqref{eq:l2-scale}; compact support gives the Lipschitz estimate
\eqref{eq:lipschitz-kernel}, permits the finite-scale transference argument
in \cref{thm:ergodic-maximal}, and is a hypothesis of Rubio de Francia's
Euclidean maximal theorem used in \cref{thm:rubio}.  No regularity of the
support is used.
\end{remark}

\section*{Statements and Declarations}

\paragraph{\textbf{Funding.}}
This work was supported by the Swedish Research Council under grant
VR 11253322.

\paragraph{\textbf{Competing interests.}}
The author has no relevant financial or non-financial interests to disclose.

\paragraph{\textbf{Data availability.}}
No datasets were generated or analysed during the current study.

\paragraph{\textbf{AI-assisted preparation.}}
During preparation of the manuscript, the author used ChatGPT (OpenAI) for
editorial assistance, consistency checking, and mathematical cross-checking,
and Aristotle (Harmonic) and Astra for automated manuscript and proof auditing.  All
AI-assisted output was reviewed and verified by the author, who takes full
responsibility for the mathematical content, references, and final text.


\end{document}